\IfFileExists{data_prediction_experiment/y_true_at_No_learning.csv}{}{%

}
\PassOptionsToPackage{dvipsnames}{xcolor}
\documentclass[final]{l4dc2026}
\usepackage{acro}
\usepackage{mathtools}
\usepackage{float}
\usepackage{amsfonts}
\usepackage{physics}
\usepackage{lilyglyphs}
\usepackage{tcolorbox}
\usepackage{epigraph} 
\usepackage{setspace}
\usepackage{matlab-prettifier}
\usepackage{textcomp}
\usepackage{mathrsfs} 
\usepackage{booktabs}
\usepackage{bm}
\usepackage{enumitem}
\usepackage{pgfplots}
\usepackage{pgfplotstable} 
\usepackage{booktabs}
\usepackage{placeins}    
\usepackage{float}       
\usepackage{wrapfig} 
\usepackage[skip=10pt]{caption}
\makeatletter

\let\For\relax
\let\EndFor\relax
\let\State\relax
\let\Require\relax
\makeatother

\usepackage{algorithm}
\usepackage[noend]{algpseudocode}  
\algrenewcommand\algorithmicindent{1.0em} 
\usepackage[framemethod=default]{mdframed}
\usepgfplotslibrary{groupplots}   
\usepgfplotslibrary{dateplot}
\pgfplotsset{compat=1.18}
\usepgfplotslibrary{fillbetween}
\usetikzlibrary{calc} 

\newcommand{\Flag}{\mathrm{Flag}}

\newcommand{\im}{\mathrm{im~}}

\newcommand{\ini}{\mathrm{ini}}

\renewcommand{\grad}{\mathrm{grad}}

\newcommand{\rmd}{\mathrm{d}}
\newcommand{\rmp}{\mathrm{p}}
\newcommand{\rmf}{\mathrm{f}}

\newcommand{\Tf}{T_\mathrm{f}}
\newcommand{\Tini}{T_\mathrm{ini}}
\newcommand{\Tsim}{T_\mathrm{sim}}
\newcommand{\Tswitch}{T_\text{switch}}

\newcommand{\Yp}{Y_{\mathrm{p}}}
\newcommand{\Yf}{Y_{\mathrm{f}}}
\newcommand{\Vp}{V_{\mathrm{p}}}
\newcommand{\Vf}{V_{\mathrm{f}}}

\newcommand{\ga}{\gamma}

\newcommand{\R}{\mathbb{R}}

\newcommand{\Zp}{\mathbb{Z}_{\geq 0}}

\newcommand{\bmU}{\mathbf{U}}
\newcommand{\bmH}{\mathbf{H}}
\newcommand{\bmV}{\mathbf{V}}

\newcommand{\bmA}{\mathbf{A}}
\newcommand{\bmB}{\mathbf{B}}

\newcommand{\calA}{\mathcal{A}}

\newcommand{\calN}{\mathcal{N}}
\newcommand{\calO}{\mathcal{O}}

\newcommand{\calT}{\mathcal{T}} 

\renewcommand{\hat}{\widehat}

\definecolor{WiscRed}{RGB}{180, 42, 26}

\definecolor{PeachRed}{RGB}{231, 143, 129}

\newcommand{\Span}{\mathrm{span}}

\DeclareMathOperator*{\diag}{diag}  

\DeclareMathOperator*{\argmin}{arg\,min}

\newcommand{\St}{\mathrm{St}}   
\newcommand{\Gr}{\mathrm{Gr}}   

\newcommand{\Exp}{\mathrm{Exp}}

\definecolor{bluecomment}{RGB}{18, 48, 174}
\definecolor{redcomment}{RGB}{217, 22, 86}
\definecolor{eqnColor}{RGB}{199, 62, 58}
\definecolor{boxColor}{RGB}{210, 224, 243}

\definecolor{myteal}{RGB}{42,160,135}   
\definecolor{myorange}{RGB}{236,112,72} 
\definecolor{mypurple}{RGB}{115,140,200}
\definecolor{mygreen}{RGB}{120,190,60}  
\definecolor{myblue}{RGB}{90,140,210}     
\definecolor{myred}{RGB}{232, 99, 99}     
\definecolor{mycyan}{RGB}{103, 196, 181}     
\definecolor{myviolet}{RGB}{170, 149, 214}

\definecolor{PeachRed}{RGB}{231, 143, 129}

\definecolor{coolgray}{RGB}{180,180,180}      
\definecolor{mutedorange}{RGB}{230,160,90}    
\definecolor{dustyrose}{RGB}{219, 64, 40}    
\definecolor{deepblue}{RGB}{90,140,210}     
\definecolor{slategreen}{RGB}{120,190,60}    
\definecolor{seateal}{RGB}{95,175,190}        
\definecolor{sandyellow}{RGB}{240, 88, 122}    
\definecolor{plumviolet}{RGB}{165,125,190}    

\DeclareAcronym{LTI}{short=LTI,long=Linear Time-Invariant}
\DeclareAcronym{LTV}{short=LTV,long=Linear Time-Varying}
\DeclareAcronym{DeePC}{short=DeePC,long=Data-enabled Predictive Control}
\DeclareAcronym{SVD}{short=SVD,long=singular value decomposition}
\DeclareAcronym{PCA}{short=PCA,long=principal component analysis}
\DeclareAcronym{FRONT}{short=FRONT,long=Flag Recursive Online Tracking}
\DeclareAcronym{GREAT}{short=GREAT,long=Grassmannian Recursive Algorithm for Tracking}
\DeclareAcronym{RGD}{short=RGD,long=Riemannian gradient descent}
\DeclareAcronym{SISO}{short=SISO,long=single-input single-output}
\DeclareAcronym{ARX}{short=ARX,long=AutoRegressive model with eXogenous inputs}
\DeclareAcronym{FROSL}{short=FROSL,long=Flag Recursive Online Subspace Learning}
\DeclareAcronym{PAST}{short=PAST,long=projection approximation subspace tracking}
\DeclareAcronym{NSR}{short=NSR,long=noise-to-signal ratio}

\title{Online Subspace Learning on Flag Manifolds for System Identification}
\usepackage{times}

\author{%
 \Name{Dian Jin} \Email{djin38@wisc.edu}\\
 \addr University of Wisconsin-Madison
 \AND
 \Name{Jeremy Coulson} \Email{jeremy.coulson@wisc.edu}\\
 \addr University of Wisconsin-Madison%
}

\begin{document}
\maketitle
\begin{abstract}%
Data-driven control methods based on subspace representations are powerful but are often limited to linear time-invariant systems where the model order is known. A key challenge is developing online data-driven control algorithms for time-varying systems, especially when the system's complexity is unknown or changes over time. To address this, we propose a novel online subspace learning framework that operates on flag manifolds. Our algorithm leverages streaming data to recursively track an ensemble of nested subspaces, allowing it to adapt to varying system dimensions without prior knowledge of the true model order. We show that our algorithm is a generalization of the Grassmannian Recursive Algorithm for Tracking. The learned subspace models are then integrated into a data-driven simulation framework to perform prediction for unknown dynamical systems. The effectiveness of this approach is demonstrated through a case study where the proposed adaptive predictor successfully handles abrupt changes in system dynamics and consistently achieves competitive performance against several baselines.
\end{abstract}


\section{Introduction}
Subspace representations are becoming increasingly important in the field of systems and control, especially with the recent revival of behavioral systems theory \citep{willems1986time} in data-driven control. A key concept in this area is Willems' Fundamental Lemma, which provides a direct data-driven non-parametric \emph{subspace} representation for finite-length trajectories of linear time-invariant systems. This has spurred the development of effective subspace identification techniques \citep{favoreel1999spc, van2012subspace}, novel data-driven simulation and control methods \citep{markovsky2008data, coulson2021distributionally, de2019formulas, berberich2020data, verhoek2021data, dorfler2022bridging}, 
which have been shown to perform remarkably well across a wide range of application domains~\citep{ elokda2021data, huang2021quadratic,  huang2021decentralized, kerkhof2021predictive, ames2022quad}.
See \citep{markovsky2021behavioral} and references therein. 

Despite these advances, many existing identification and control methods are limited to \ac{LTI} systems. A significant open challenge is the design of online adaptive approaches and recursive algorithms that can handle systems that change over time. While some subspace tracking methods have been developed for online state-space identification \citep{verhaegen2007filtering, van1994n4sid, oku2002recursive, zhang2019online, sasfi2025GREAT}, they typically require the system order to be known and constant. In many real-world scenarios, however, the true system order is unknown and may vary, for instance, due to component failures. 
In this work, we propose a novel framework for adaptive, data-driven identification of unknown, time-varying dynamical systems based on subspace learning techniques.

Subspace learning (sometimes termed subspace tracking) has long been studied in signal processing \citep{yang1995projectionPAST, delmas2010subspace, balzano2010GROUSE, balzano2015local, he2012incremental}.  
Broadly speaking, subspace learning algorithms can be classified into algebraic and geometric methods \citep{balzano2018streamingPCA}. The geometric methods optimize a certain loss function via gradient descent on a matrix manifold such as the Grassmannian, a set consisting of all subspaces of a particular dimension. Such methods have a solid theoretical background rooted in optimization on manifolds \citep{edelman1998geometry, absil2004riemannian, boumal2023introduction}. 
Running these algorithms on Grassmannians requires choosing the dimensionality \emph{a priori}. This poses a significant challenge in practical applications where the true underlying dimension is unknown or may change over time. 
\paragraph{Contribution.}
We propose a novel framework for online system identification and data-driven simulation for unknown time-varying dynamical systems. The proposed framework consists of the \ac{FRONT} algorithm on flag manifolds that leverages streaming data and yields a nested hierarchy of subspaces.
To the best of our knowledge, this is the first online algorithm that optimizes over flag manifolds with streaming data.
We formally show that \ac{GREAT}~\citep{sasfi2025GREAT} is a special case of our algorithm when the true subspace dimension is known. When the dimension is unknown, we propose a dimensionality ensemble strategy inspired by ensemble learning \citep{dietterich2000ensemble}.  
Finally, the \ac{FRONT} algorithm is leveraged as a subspace model in data-driven simulation \citep{markovsky2008data} to perform trajectory prediction for a time and complexity-varying unknown system.
The framework is showcased on an \ac{ARX} system simulation and is shown to outperform non-adaptive subspace prediction \citep{favoreel1999spc} and other adaptive benchmarks. 

The rest of the paper is organized as follows. Section~\ref{sec:preliminaries} contains preliminaries on subspace geometry. 
We formulate the problem in Section \ref{sec:problem-setup}.
Section~\ref{sec:flag_tracking_algorithm} introduces the online subspace learning algorithm.  
Section \ref{sec:numrical-study} contains numerical studies.
Section~\ref{sec:conclusion} concludes the paper.
\section{Preliminaries}\label{sec:preliminaries}
We use $\Zp$ to denote non-negative integers.
Given $i,j\in\Zp$, with $i\leq j$, we use $[i, j]$ to denote the discrete interval $[i, j]\cap \Zp$. 
 We use $I_{p \times m}$ to denote the $p \times m$ matrix whose entries are $1$ on the main diagonal (i.e., $(i, i)$-entry for $i=1,\dots,\min(p, m)$), and all other entries are $0$. 
 The diagonal matrix with diagonal elements $\ga_1,\cdots,\ga_d$ is denoted $\diag(\ga_1,\cdots,\ga_d)$. 
 Given a function  $v\colon \Zp\to\R^q$ and $i,j\in\Zp$, with $i\leq j$, we define $v_{[i,j]}  =  (v(i), v(i+1), \dots, v(j))  \in \R^{q(j-i+1)}$. We use $\bmU$ to denote a subspace of $\R^p$ and $U$ to denote a matrix such that $\im U = \bmU$. The orthogonal projector onto $\bmU$ is denoted by $\Pi_{\bmU}$, and is used interchangeably with $\Pi_{U} = UU^\top$, where $U$ is an orthonormal matrix whose columns span $\bmU$. 
 The Moore–Penrose inverse of a matrix $A$ is denoted $A^\dag$. 
 The Gaussian distribution with mean $0$ and covariance $\Sigma \in \R^{p \times p}$ is denoted $\calN(0, \Sigma)$.
 
\subsection{Subspace geometry}
A Grassmannian is the set of all subspaces in $\R^p$ with dimension $d$ denoted by $\Gr(p, d)\coloneqq \{\bmU $ a linear subspace of $\R^p: \dim \bmU = d \}.$
Flag manifolds generalize the Grassmannians by considering sequences of nested subspaces of increasing dimension.
Let $p \geq 2$ be an integer and $q_{1:d} \coloneqq (q_1, \dots, q_d)$ with $0 < q_1 < \cdots < q_d < p.$
A flag of signature $(p, q_{1:d})$ is a nested sequence of linear subspaces $\{\bmU_i\}_{i=1}^{d}$ of $\R^p$ satisfying
$\{0\} \subset \bmU_1 \subset \cdots \subset \bmU_d \subset \R^p$
with $\dim \bmU_j = q_j, j=1,\dots,d$. We denote this flag by $\bmU_{1:d}$. 
The set of all such flags forms a smooth manifold $\Flag(p, q_{1:d})$ \citep{ye2022optimization}.
We also denote $q_0 = 0, q = q_d, q_{d+1} = p$. 
We can represent flags as points on the Stiefel manifold $\St(p, q) = \{U \in \R^{p \times q}: U^\top U = I\}$ consisting of orthonormal matrices. 
Define for each $k =1,\dots,d$, $U_k \in \St(p, q_{k}-q_{k-1})$ such that $[U_1  \cdots  U_k]$ is an orthonormal basis of $\bmU_k$. We must note that under these notations, $\bmU_k$ is not the subspace spanned by $U_k$.
Then $U_{1:d} \coloneqq [U_1  \cdots  U_d]$ is a representative of $\bmU_{1:d}$, called Stiefel representative. Throughout this paper we will abuse the notation to write $U_{1:d} \in \Flag(p, q_{1:d})$. 
If the signature is $(p, q_1)$, then $\Flag(p, (q_1)) = \Gr(p, q_1)$.
The chordal distance of two subspaces $\bmA \in \Flag(p, (m)), \bmB \in \Flag(p, (n))$ is defined by $\rmd(\bmA, \bmB) = (\sum_{i=1}^{\min(m,n)} \sin^2 \theta_i)^{1/2}$ \citep[Theorem 12]{ye2016schubert}, where $\theta_i$ is the $i^{\text{th}}$ principal angle \citep{golub2013matrix}
between $\bmA$ and $\bmB$. 
We present a simple example to illustrate the difference between flags and subspaces.
\begin{example}
Consider the canonical basis $e_1=(1, 0, 0, 0), e_2 = (0, 1, 0, 0), e_3 = (0, 0, 1, 0), e_4 = (0, 0, 0, 1)$ of $\R^4$. Let $\bmU_1 = \Span \{e_1\}, \bmU_2 = \Span \{e_1, e_2\}, \bmU_3 = \Span \{e_1, e_2, e_3\}$, then $\bmU_{1:3}\in \Flag(4, (1, 2, 3))$, and a Stiefel representation of this flag is the matrix $U = \begin{bmatrix} e_1 & e_2 & e_3 \end{bmatrix}$. 
Note that even if two bases span the same subspace, the hierarchical structures of their flags may be different. Let $V = \begin{bmatrix} e_1 & e_3 & e_2 \end{bmatrix}$, then clearly $\im V = \im U$. However, $\bmV_1 = \Span\{e_1\}, \bmV_2 = \Span\{e_1, e_3\}, \bmV_3 = \Span\{e_1, e_2, e_3\}$, which implies $\bmV_2 \neq \bmU_2$, hence $\bmV_{1:d}\neq \bmU_{1:d}$. 
\end{example}
\section{Problem Setup \& Motivation}\label{sec:problem-setup}
The goal of this paper is to estimate an unknown, time-varying subspace $\bmU^t$ from streaming data. More precisely, for each $t\in\Zp$, let
$\bmU^t \in \bigcup_{j=1}^{d} \Gr(p, q_j),$
where $0<q_1<\cdots<q_d<p$. We assume that we have access to an upper bound $q_d$ on the dimension of all subspaces $\bmU^t$. Note that the upper bound $q_d$ may not be tight in the sense that there may be no time instance $t\in\Zp$ where $\dim \bmU^t = q_d$. We assume that at each $t\in\Zp$, we measure the data $w_t=v_t+\eta_t,$
where $v_t\in \bmU^t$, and $\eta_t$ is noise.
We wish to obtain estimates $\hat{\bmU}^t$ such that
$d(\hat{\bmU}^t,\bmU^t)\leq\epsilon_{\textup{tol}}$
for some user-defined tolerance $\epsilon_{\textup{tol}}>0$, i.e., we seek estimates that are ``close'' to the true subspaces. 
We refer to this problem as \emph{learning} $\bmU^t$ or \emph{tracking} $\bmU^t$.
In the case when $d=1$ and $q_d$ is known, there are efficient algorithms~\citep{yang1995projectionPAST, balzano2010GROUSE, sasfi2025GREAT} for robustly learning $\bmU^t$. However, in many practical settings the true dimension of the underlying subspaces is \emph{unknown} and may even vary over time (e.g., online dynamic mode decomposition for time-varying dynamical systems~\citep{zhang2019online}, direction of arrival analysis \citep{rabideau1996fast}). 
Solving this problem is critical for online system identification and control: a fixed, mis-specified subspace dimension yields biased estimates, slow adaptation, and brittle decisions \citep{ljung1998system} (see Figure \ref{fig:y_prediction_two}). Learning a nested hierarchy of subspaces (flag) online gives rank-on-demand models that track changing complexity, enabling efficient prediction, reliable change-point detection, without prior knowledge of the true order \citep{szwagier2025nestedsubspacelearningflags}.
\paragraph{Motivating Example: Online System Identification.}\label{sec:motivation}
Consider the \ac{LTV} dynamical system 
\begin{equation}\label{eq:LTV}
\begin{aligned}
    x(t+1) &= A_t x(t) + B_t u(t), \\
    y(t) &= C_t x(t) + D_t u(t),
\end{aligned}
\end{equation}
where $A_t\in\R^{n\times n}$, $B_t\in\R^{n\times m}$, $C_t\in\R^{r\times n}$, and $D_t\in\R^{r\times m}$ are time varying matrices. Assume that the model parameters $(A_t,B_t,C_t,D_t)$ are unknown, and we only have access to input-output data of~\eqref{eq:LTV}. 
A typical problem in system identification is to use input-output data to find a model that captures the input-output behavior of~\eqref{eq:LTV}. The main challenge in this setting is that the model parameters are \emph{time-varying}. Furthermore, we do not have access to state measurements. Hence, the underlying state-space dimension (system complexity) is unknown.
Let $L \geq 1,$ for each $t  \geq L-1$, the $L$-length input-output trajectories of~\eqref{eq:LTV} on the interval $[t-L+1,t]$ are given by
\begin{equation*}
v_t := \begin{bmatrix}
    u_{[t-L+1,t]} \\ y_{[t-L+1,t]}
\end{bmatrix} = \underbrace{\begin{bmatrix}
    0 & I_{mL} \\
    \mathcal{O}_{[t-L+1,t]} & \mathcal{T}_{[t-L+1,t]}
\end{bmatrix}}_{=: U^t}
\begin{bmatrix}
    x(t-L+1) \\ u_{[t-L+1,t]}
\end{bmatrix}, \text{ where}
\end{equation*}
\[
\resizebox{\columnwidth}{!}{$\calO_{[t-L+1,t]} = \begin{bmatrix}
    C_{t-L+1} \\
    C_{t-L+2} A_{t-L+1} \\
    \vdots \\
    C_{t} A_{t-1}\cdots A_{t-L+1}
\end{bmatrix},
\calT_{[t-L+1,t]} =\begin{bmatrix}
    D_{t-L+1} & 0 & 0 & \cdots & 0 \\
    C_{t-L+2} B_{t-L+1} & D_{t-L+2} & 0 & \cdots & 0 \\
    C_{t-L+3} A_{t-L+2} B_{t-L+1} & C_{t-L+3} B_{t-L+2} & D_{t-L+3} & \cdots & 0 \\
    \vdots & \vdots & \vdots & \ddots & \vdots \\
    C_{t} A_{t-1} \cdots A_{t-L+2}B_{t-L+1} & \cdots & \cdots & \cdots & D_{t}
\end{bmatrix}.$}
\]
Instead of learning the state-space matrices $(A_t, B_t, C_t, D_t)$ through (recursive) state-space system identification methods~\citep{verhaegen1992subspace, van1994n4sid, oku2002recursive}, we directly learn the subspaces $\bmU^t \in \bigcup_{j=1}^d \Gr(p, q_j)$ based on streaming data $w_t = v_t + \eta_t$, where $\eta_t$ is measurement noise. Note that $\dim \bmU^t$ is unknown and may vary depending on the rank of $\calO_{[t-L+1,t]}$.
The subspace $\bmU^t$ is known in behavioral system theory \citep{willems1986time, markovsky2021behavioral} as the \textit{restricted behavior} \citep{markovsky2022identifiability} and serves as the backbone of several powerful nonparametric control methods~\citep{coulson2019data, de2019formulas, markovsky2021behavioral, markovsky2008data}.
This online learning problem is studied in \citep{sasfi2025GREAT} under the assumption that $\dim \bmU^t$ is constant, which may not hold in practical settings in which the complexity of the system varies over time (see Section \ref{subsec:geodesic_tracking}). 
In what follows, we propose a method for learning $\bmU^t$ without assuming $\dim \bmU^t$ is constant, and use the learned flags for a data-driven prediction task (see Section \ref{sec:ARX-system}). 
While we do not address control here, the learned subspaces can be incorporated into data-driven control frameworks, which is a direction for future work.
\section{Online subspace learning on flag manifolds}\label{sec:flag_tracking_algorithm}
Inspired by recent advances in subspace learning for online system identification \citep{sasfi2025GREAT}, we propose a method for solving the learning problem in Section \ref{sec:problem-setup}, which consists of estimating the evolving subspace $\bmU^t$ using optimization on the flag manifold. Each element of the flag manifold is an ensemble of subspaces that can be used for downstream tasks such as data-driven prediction. Thus, our algorithm contributes to a data-to-prediction pipeline that continuously adapts to streaming data and the changing complexity of the system (see Section \ref{sec:ARX-system}).
\subsection{Gradient descent on flag manifolds with streaming data}\label{subsec:why_flag}
Let $\{\bmU^t\}_{t \in \Zp}$ be a sequence of subspaces with each $\bmU^t \in \bigcup_{j=1}^{d} \Gr(p, q_j)$. At each time $t \in \Zp$ we observe a data sample $w_t = v_t + \eta_t$, where $v_t \in \bmU^t$ and $\eta_t$ is noise.
Let $T \geq 1$. For each $t \geq T-1$, let $W_t = [w_{t-T+1}~\cdots~w_t] \in \R^{p \times T}$ be the window of the most recent $T$ data samples. 
For each $t \in \Zp$ and $j\in\{1,\dots,d\}$, we define the cost function $g_{W_t}^j \colon \Gr(p,q_j)\to\R$ by
\begin{equation}\label{eq:grass_cost_function}
g_{W_t}^j(\bmU) = \norm{W_{t} - \Pi_{\bmU}W_{t}}_F^2,
\end{equation}
where $\| \cdot \|_F$ is the Frobenius norm. 
Minimizing~\eqref{eq:grass_cost_function} over $\bmU \in \Gr(p, q_j)$ is known as \ac{PCA}~\citep[Section 2.4]{boumal2023introduction} and entails finding the $q_j$-dimensional subspace that best fits the data in the sense of minimizing the projection error. 
While alternative norms (e.g., $\ell_1$) could provide improved robustness under non-Gaussian noise \citep{candes2011robust}, our focus in this work is on the geometric and algorithmic aspects of flag-based subspace tracking. 
It is well-known that \ac{PCA} admits a closed-form solution~\citep{eckart1936approximation} given by the top $q_j$ singular vectors (principal vectors) of $W_t$. 
One approach to solving the proposed problem would be running $d$ separate optimizations on each $\Gr(p, q_k)$ and select a suitable dimension, but this is computationally expensive and would yield optimal subspaces $\bmV_j^t := \argmin_{\bmU \in \Gr(p, q_j)} g_{W_t}^j(\bmU)$ that are not nested for $j=1, \dots,d$.
The lack of nestedness has been shown \citep{szwagier2025nestedsubspacelearningflags} to cause poor performance when using subspaces for tasks such as robust subspace recovery. 
Instead, we propose optimizing over $\Flag(p,q_{1:d})$.
For each $t\in\Zp$ we define
\begin{equation}\label{eq:flag_cost_function}
    f_{W_t}(\bmU_{1:d}) = \norm{W_{t} - \frac{1}{d}\sum_{i=1}^{d}\Pi_{\bmU_i}W_{t}}_F^2
\end{equation}
for all $\bmU_{1:d} \in \Flag(p, q_{1:d})$. 
One key difference of optimizing~\eqref{eq:flag_cost_function} instead of~\eqref{eq:grass_cost_function} is that we obtain a nested sequence of subspaces which can be leveraged for downstream tasks (e.g., data-driven prediction) via ensembling methods~\citep{dietterich2000ensemble} (see Section~\ref{sec:ARX-system}). 
The operator $\frac{1}{d}\sum_{i=1}^{d}\Pi_{\bmU_i}$ demonstrates an ensemble of projections onto subspaces of different dimensions, which is called the ``flag trick'' in \citep[Definition 5]{szwagier2025nestedsubspacelearningflags}. Similar to~\eqref{eq:grass_cost_function}, minimizers of \eqref{eq:flag_cost_function} can be computed in closed form by computing the top-$q_d$ principal vectors of $W_t$~\citep[Theorem 4]{szwagier2025nestedsubspacelearningflags}. 
However, 
recomputing the full sample covariance $W_t W_t^\top$ and its leading eigenspaces at every time $t$ can be computationally expensive, 
motivating recursive subspace tracking methods \citep{balzano2018streamingPCA}. 
Based on results in optimization on flag manifolds \citep{ye2022optimization, zhu2024practical}, we perform \ac{RGD} \citep[Section 4.3]{boumal2023introduction} for \eqref{eq:flag_cost_function}. 
The goal of \ac{RGD} is to seek a minimizer of the function \eqref{eq:flag_cost_function} by iterative updates. In particular, for each new data sample $w_t$, we perform $K>0$ steps along the gradient of~\eqref{eq:flag_cost_function}. Suppose that the initial point is $\hat{\bmU}^{T-1}_{1:d}$ (we use the superscript $T-1$ since we have to collect $T$ data points to form the first window $W_{T-1}$). For each fixed $t \geq T-1$, 
we update
\begin{equation}\label{eq:flag_update}
\begin{aligned}
\hat{\bmU}^{t, k+1}_{1:d} = \Exp_{\hat{\bmU}^{t, k}_{1:d}}
\left(-\alpha_t \grad~f_{W_t}(\hat{\bmU}^{t, k}_{1:d})\right),
\end{aligned}
\end{equation}
where $k=0,\dots,K-1$, $\alpha_t > 0$ is the step size, and we denote 
$\hat{\bmU}^{T-1,0}_{1:d} = \hat{\bmU}^{T-1}_{1:d}$, 
$\hat{\bmU}^{t+1}_{1:d} \coloneqq \hat{\bmU}^{t, K}_{1:d}$,
$\hat{\bmU}^{t+1,0}_{1:d} = \hat{\bmU}^{t+1}_{1:d}$, for all $t \geq T-1$. 
The Riemannian gradient $\grad~f_{W_t}$ is a vector field on the manifold, and $\Exp_{\hat{\bmU}^{t, k}_{1:d}}$ is the exponential map, which maps the tangent vector $-\alpha_t \grad~f_{W_t}(\hat{\bmU}^{t, k}_{1:d})$ back to the manifold. 
Riemannian gradient and exponential map can be computed using only matrix algebra \citep{ye2022optimization}.

\subsection{Learning time- and dimension-varying subspaces}\label{subsec:flag_tracking_algorithm}
Having established the necessary foundations,
we now present the \textbf{F}lag \textbf{R}ecursive \textbf{ON}line \textbf{T}racking (FRONT) algorithm given by Algorithm \ref{alg:flag}.
\begin{figure}[tbp]
\centering
\begin{minipage}{0.48\linewidth}
\begin{algorithm}[H]
\small
\caption{FRONT}\label{alg:flag}
\begin{algorithmic}[1]
  \Require Window length $T \geq 1$; initial estimate $\hat{\bmU}^{T-1}_{1:d}$; samples $\{w_t\}_{t \in \Zp}$; step sizes $\alpha_t$
  \For{$t = T-1,T,\dots$}
    \State Construct $W_t=[w_{t-T+1}\cdots w_t]$
    \Statex \quad Initialize gradient descent: $\hat{\bmU}^{t,0}_{1:d} = \hat{\bmU}^{t}_{1:d}$
    \For{$k = 0,1,\dots,K-1$}
      \Statex \qquad  Gradient step: 
      \Statex \qquad  $\hat{\bmU}^{t,k+1}_{1:d} =
        \Exp_{\hat{\bmU}^{t,k}_{1:d}}\!\left(-\alpha_t \grad f_{W_t}(\hat{\bmU}^{t,k}_{1:d})\right)$

    \EndFor
    \State \textbf{Update:} $\hat{\bmU}^{t+1}_{1:d} = \hat{\bmU}^{t,K}_{1:d}$
  \EndFor
\end{algorithmic}
\end{algorithm}
\end{minipage}
\hfill
\begin{minipage}{0.48\linewidth}

\begin{algorithm}[H]
\small
\caption{GREAT}\label{alg:GREAT}
\begin{algorithmic}[1]
  \Require Window length $T \geq 1$; initial estimate $\hat{\bmU}^{T-1}$; samples $\{w_t\}_{t \in \Zp}$; step sizes $\alpha_t$
  \For{$t = T-1,T,\dots$}
    \State Construct $W_t=[w_{t-T+1}\cdots w_t]$
    \Statex \quad Initialize gradient descent: $\hat{\bmU}^{t,0} = \hat{\bmU}^{t}$
    \For{$k = 0,1,\dots,K-1$}
      \Statex \qquad Gradient step:
      \Statex \qquad $\hat{\bmU}^{t,k+1} = \Exp_{\hat{\bmU}^{t,k}}\!\left(-\alpha_t\grad f_{W_t}(\hat{\bmU}^{t,k})\right)$
    \EndFor
    \State \textbf{Update:} $\hat{\bmU}^{t+1} = \hat{\bmU}^{t,K}$
  \EndFor
\end{algorithmic}
\end{algorithm}
\end{minipage}
\end{figure}
The output of \ac{FRONT} is a sequence of flags, rather than selecting a single subspace from $\bmU_{1:d}^t$.  
We aggregate information across subspaces of different dimensions, with their usage depending on the specific task.
Since a flag is represented by an orthonormal matrix in $\R^{p \times q_d}$, the per-time-step complexity of \ac{FRONT} is $O(K p^2 q_d)$, by the same covariance-based complexity argument as in \citep[Remark 5]{sasfi2025GREAT}.
FRONT is similar in structure to \ac{GREAT} \citep{sasfi2025GREAT}, but generalizes \ac{GREAT} by learning an ensemble of models with different complexity while not assuming that the true data-generating subspaces have fixed dimension.
Moreover, when $d=1$ so that $\Flag(p, (q_1)) = \Gr(p, q_1)$, \ac{FRONT} is equivalent to \ac{GREAT}. To prove this, we make a definition on algorithm equivalence. 
\begin{definition}\label{def:algorithm_equivalence}
    Let $M$ be a Riemannian manifold.
    We define an algorithm $\calA$ on $M$ as an update mapping: $\calA: M \to M$ by $U^{t+1} = \calA(U^t).$
\end{definition}
\begin{definition}
     Let $\calA_1, \calA_2$ be two algorithms on $M$ such that $U^{t+1}=\calA_1(U^t)$ and $V^{t+1}=\calA_2(V^t)$, respectively. 
    We say that $\calA_1$ is equivalent to $\calA_2$ if 
    $U^0=V^0 \text{ implies } U^t = V^t \text{ for all } t \geq 1. $
\end{definition}

\begin{proposition}
Let $\calA_1: \Flag(p, (q)) \to \Flag(p, (q))$ be defined by (\ref{eq:flag_update}), and $\calA_2:\Gr(p,q) \to \Gr(p,q)$ be defined by the update rule in \ac{GREAT}. 
Then, $\calA_1$ is equivalent to $\calA_2$.
\end{proposition}
\begin{proof}
    By definition, $\Flag(p, (q)) = \Gr(p, q)$, so $\calA_1$ and $\calA_2$ are defined on the same manifold. 
    Let $\bmU^0$ be the initial point of $\calA_1$ and $\bmV^0$ be the initial point of $\calA_2$. Suppose $\bmU^0 = \bmV^0$.
    The cost function (\ref{eq:flag_cost_function}) becomes 
    $f_{W_t}(\hat{\bmU})=\|W_t - \Pi_{\hat{\bmU}}W_t\|_F^2$, which coincides with $g_{W_t}(\hat{\bmU})$. For each $t = T-1, T, \dots$ and $k \in \{0, \dots,K-1\}$,
    let $\bmU^{t,k}$ and $\bmV^{t,k}$ be the points obtained by 
    the gradient step in \ac{FRONT} and \ac{GREAT}, respectively.
    Since $\Flag(p, (q))$ and $\Gr(p, q)$ are the same manifold, their intrinsic geometric structures coincide. Thus, the exponential map and the Riemannian gradients are identical. This can also be verified by the formulas in \citep{ye2022optimization} for $\Flag(p,(q))$ with those in \citep{edelman1998geometry} for $\Gr(p,q)$.
    Therefore, we have $\hat{\bmU}^{t, k} =  \hat{\bmV}^{t, k}$ for all $t \in \{T-1,T,\dots\}$ and $k \in \{0,\cdots,K-1\}$, hence $\calA_1$ and $\calA_2$ are equivalent.
\end{proof}
    Under assumptions on  quantitative persistency of excitation of the data and bounded subspace variation~\citep[Assumptions 1-4]{sasfi2025GREAT},  all convergence guarantees and finite-sample analysis of \citep[Theorem 1]{sasfi2025GREAT} directly extend to \ac{FRONT} on $\Flag(p, (q))$. 
    In particular, \citep {sasfi2025GREAT} derives finite-sample upper bounds on the distance $d(\hat{\bmU}^t,\bmU^t)$ in terms of subspace change rate, and singular values of the data window $W_t$. These results provide guaranteed convergence rates to the true underlying subspace and therefore give a quantitative method for determining conditions on data needed to satisfy $d(\hat{\bmU}^t,\bmU^t) < \epsilon_{\text{tol}}$ for a given $\epsilon_{\text{tol}}$. 
    Beyond this trivial signature, we do not claim convergence guarantees for \ac{FRONT} in this paper. 
    In the future work, we wish to study finite-sample convergence guarantees for \ac{FRONT} in the more general case when the flag does not have a trivial signature $(p, (q))$.
    The nontrivial flag-manifold setting is supported by the following numerical case studies.
\section{Numerical case studies}\label{sec:numrical-study}
Throughout all numerical examples\footnote{The code is available at \url{https://github.com/DianJin-Frederick/FRONT-experiments}.}, when implementing \ac{FRONT}, we use the Python class \texttt{Flag} developed in \citep{szwagier2025nestedsubspacelearningflags} and use the \texttt{SteepestDescent} optimizer from Pymanopt \citep{pymanopt}, which uses backtracking line-search~\citep{boumal2023introduction} to specify the step sizes $\alpha_t$ used for the gradient steps in \ac{FRONT}. The number of gradient steps is $K=5$. 
\subsection{Geodesic tracking}\label{subsec:geodesic_tracking}
We generate a sequence of time-varying ground-truth subspaces
on $\Flag(10, (1,\dots,5))$. The initial flag $\bmU^0$ has Stiefel representation
$
U^{0} = I_{10\times 5}.
$
For each $t \in \Zp$ we update $\bmU^{t+1} = \Exp_{\bmU^t}(\alpha \bmH^t),$
where $\bmH^t$ is a random tangent vector \citep[Section 8.3]{boumal2023introduction} at $\bmU^t$
and $\alpha = 5 \cdot 10^{-5}$.
At time \(T_{\text{switch}}=100\), we increase the subspace dimension by setting
$
U^{\Tswitch} = I_{10 \times 6}
$
and generate subsequent flags in the same way on $\Flag(10, (1,\dots,6))$. 
At each time $t$ we observe a data point $w_t = U^t a_t + \eta_t$, where $a_t \sim \calN(0, I)$ and $\eta_t \sim \calN(0, 10^{-4}I)$ is noise.
For each $t \geq T-1$, the data window is constructed as $W_t = [w_{t-T+1} \dots w_t]$.
The subspace estimate is updated online by running \ac{FRONT} on $\Flag(10, (5, 6))$.
We set the initial estimate $\hat{\bmU}^T$ to be random and run with data window size $T \in \{1, 20, 50\}$. For each $T \in \{1, 20, 50\}$, we run 100 experiments across different random data sets $\{w_t\}_{t=1}^{200}$. 

\begin{wrapfigure}[8]{r}{0.45\linewidth} 
\vspace{-10pt} 
\centering
\begin{tikzpicture}
\begin{axis}[
    xmin=0, xmax=200, xtick={0,20, 50,100,150,200},
    ymin=0, ymax=1.30, 
    width=1.0\linewidth,
    height=0.5\linewidth,
    tick label style={font=\scriptsize},
    label style={font=\scriptsize},
    legend style={
        draw=gray!0,
        font=\scriptsize,
        at={(0.5,1.0)}, anchor=south,
        legend columns=3, 
        /tikz/every even column/.append style={column sep=6pt}
    },
    grid=both,
    cycle list={
        {myteal,   mark=o},
        {myorange, dashed, mark=square*},
        {mypurple, dashed, mark=triangle*}
    }
  ]
    \addplot[thick, color=myteal, line width=1.0pt] table[
      col sep=comma,
      header=false,
      y index=0,
      x expr=\coordindex
    ]{geodesic_tracking_data/chordal_distances_1_flag_5to6.csv};
    \addlegendentry{$T=1$}
    \addplot[thick, color=myorange, line width=1.0pt] table[
      col sep=comma,
      header=false,
      y index=0,
      x expr=\coordindex+20
    ]{geodesic_tracking_data/chordal_distances_20_flag_5to6.csv};
    \addlegendentry{$T=20$}
    \addplot[thick, solid, color=mypurple, line width=1.3pt] table[
      col sep=comma,
      header=false,
      y index=0,
      x expr=\coordindex+50
    ]{geodesic_tracking_data/chordal_distances_50_flag_5to6.csv};
    \addlegendentry{$T = 50$}
  \end{axis}
\end{tikzpicture}
\captionsetup{width=0.85\linewidth}
\caption{Average chordal distance for different window sizes across 100 experiments.}
\label{fig:chordal_distance}
\end{wrapfigure}
Figure \ref{fig:chordal_distance} shows how $d(\hat{\bmU}^{t}, \bmU^t)$ evolves over time $t=T-1,\dots,200$ for different window sizes $T \in \{1, 20, 50\}$. For all $T\in\{1, 20, 50\}$, the distance decreases as more data are available, indicating that the estimates gradually align with the true subspaces. At time $\Tswitch$, the distance sharply increases due to the sudden change in subspace dimension, after which \ac{FRONT} readapts and the distance decreases again. 
The plateau for $T=1$ after switching arises because a single-sample update does not provide sufficient information to identify the new subspace, whereas
larger data window sizes ($T=20, 50$) yield better estimates. 
Although this example considers an increase in subspace dimension, the problem setup in Section \ref{sec:problem-setup} does not assume increasing complexity. Decreasing subspace dimension can likewise be accommodated, since each flag contains all lower-dimensional nested subspaces.

\subsection{Online data-driven prediction for a switched system}\label{sec:ARX-system}
The data-driven prediction problem posed in \citep[Section 4]{markovsky2008data} aims to find the output sequence of an unknown dynamical system corresponding to an input sequence and a past input-output trajectory collected from the system. For \ac{LTI} systems, data-driven prediction can be performed using a matrix representation of the set of finite-length trajectories (restricted behavior)~\citep{favoreel1999spc}. Motivated by the subspace representations of dynamical systems in Section~\ref{sec:motivation}, we use \ac{FRONT} to learn the subspace of finite-length trajectories of a time- and complexity-varying system and perform data-driven prediction.
We consider the data-driven prediction task for a \ac{SISO} \ac{ARX} system, whose dynamics is given by 
\begin{equation}\label{eq:ARX}
\begin{cases}
    y_t = 0.3 y_{t-1} - 0.02 y_{t-2} + 0.6 u_{t-1} + 0.2 u_{t-2}, & t<\Tswitch, \\
    y_t = 1.5 y_{t-1} - 0.74 y_{t-2} + 0.12 y_{t-3} + 0.6 u_{t-1} + 0.2 u_{t-2} + 0.05 u_{t-3}, & t\geq \Tswitch, 
    \end{cases}
\end{equation}
The input dimension is $m=1$ and output dimension is $n=1$.
The main purpose of this example is to test adaptation to a change in system complexity, since the underlying system order increases at $\Tswitch$. This is precisely the setting that motivates the use of a flag, as it incorporates information from subspaces across all candidate dimensions.
Let $\Tsim=300$ and $\Tswitch=100$. Given a random input sequence $u_{[0,\Tsim-1]}$ and a prediction horizon $\Tf \in \Zp$, we wish to predict the corresponding future output trajectory $y_{[0,\Tsim-1]}$ of~\eqref{eq:ARX}. We assume that the measurement of output at time $t$ is corrupted by the independent noise $\eta_t$ for $t\in [0,\Tsim - 1]$. 
The \ac{NSR} at time $t \in [0, \Tsim - 1]$ is defined as $\sigma_t \coloneqq \|\eta_t\|_2 / \|y_t\|_2$. 
In our experiments, we fix this ratio to a constant $\sigma > 0$ for all $t$. The measurement noise is generated as $\eta_t \sim \mathcal{N}\!\left(0,\, (\sigma \|y_t\|_2)^2 I\right)$.
\paragraph{Offline Data Collection.}
Assume that we have access to a noisy input-output trajectory of~\eqref{eq:ARX} of length $T_\rmd=30$ denoted by $(u^\rmd_{[0,T_\rmd-1]}, y^\rmd_{[0,T_\rmd-1]} + \eta_{[0,T_\rmd]}^\rmd)$, where $\eta_t^\rmd \sim \mathcal{N}\!\left(0,\, (\sigma \|y_t^\rmd\|_2)^2 I\right)$ for $t \in [0, T_\rmd-1]$.
Let $\Tini=4$, $\Tf=4$, and $L=\Tini+\Tf=8$. 
We sort the offline input-output data into a Hankel matrix
\[
    H \coloneqq \begin{bmatrix}
    u^\rmd_{[0, L-1]} & \cdots & u^\rmd_{[T_\rmd-L, T_\rmd-1]} \\[5pt]
    y^\rmd_{[0, L-1]} + \eta_{[0,L-1]}^\rmd & \cdots & y^\rmd_{[T_\rmd-L, T_\rmd-1]} + \eta_{[T_\rmd-L, T_\rmd-1]}^\rmd 
\end{bmatrix}
\]
with the thin \ac{SVD} being $H = U \Sigma V^\top$, where $U \in \R^{(m+n)L \times r}$ is orthonormal with $r = \rank H$. Let $U = [u_1 \dots u_r]$, we set the initial point of \ac{FRONT} to be $\bmV_{1:d}^0 = \{\bmV^0_{k}\}_{k=1}^d \in \Flag((m+n)L, (q_1,\dots,q_d))$, where $\bmV^0_{k} = \Span\{u_1, \cdots, u_{q_k}\}$ for $k \in [1,d]$, and $q_d=r$. 
\paragraph{Subspace Predictor.}
Let $V \in \R^{(m+n)L \times r}$ be any orthonormal matrix with the partition
$V=[\Vp^\top ~ \Vf^\top ~ \Yp^\top ~ \Yf^\top]$,
where 
$\Vp \in \R^{m\Tini \times r}, \Vf \in \R^{m\Tf \times r}, \Yp \in \R^{n\Tini \times r}, \Yf \in \R^{n\Tf \times r}$. 
The $\Tf$-length predicted output at time $t$ is given by
\begin{equation}\label{eq:subspace-predictor}
    y^{\text{pred}}_{[t, t+\Tf-1]}(V) = 
    \begin{bmatrix}
    y^{\text{pred}}_t(V) \\ \vdots \\ 
    y^{\text{pred}}_{t+\Tf-1}(V)
    \end{bmatrix}
    = Y_\rmf 
    \begin{bmatrix}
    \Vp \\ 
    \Vf \\
    Y_\rmp 
    \end{bmatrix}^\dag
    \begin{bmatrix}
    u_{[t-T_\ini, t-1]} \\ u_{[t,t+\Tf-1]} \\ 
    y_{[t-T_\ini, t-1]} + \eta_{[t-\Tini, t-1]}
    \end{bmatrix} \in \R^{n\Tf},
\end{equation}
which is known as the orthonormal subspace predictor \citep{jin2026sensitivity}.
In order to cope with the changing dynamics of~\eqref{eq:ARX}, we use \ac{FRONT} to update the subspace (restricted behavior of \eqref{eq:ARX}: see Section \ref{sec:motivation}) based on streaming data $w_t \coloneqq (u_{[t-T+1, t]}, y_{[t-T+1, t]} + \eta_{[t-T+1, t]})$ and initial point $\bmV_{1:d}^0$, where $\eta_t \sim \calN(0, (\sigma \|y_t\|_2)^2 I)$. We then update the matrix $Y_\rmf([\Vp^\top ~\Vf^\top~ Y_\rmp^\top]^\top)^\dag$ in \eqref{eq:subspace-predictor} based on the learned flag.
Let $\hat{\bmV}_{1:d}^t = \{\hat{\bmV}_k^t\}_{k=1}^d$ be the $t^\text{th}$ flag estimate produced by \ac{FRONT}.
The Stiefel representative $\hat{V}_{1:k}^t$ of $\hat{\bmV}_k^t$ is obtained by taking the first $q_k$ columns of the Stiefel representative $\hat{V}_{1:d}^t$ of $\hat{\bmV}_{1:d}^t$. Then we partition the rank-$q_k$ matrix $\hat V_{1:k}^{t}$ as
$
\hat V_{1:k}^{t} = [V_{\rmp,k}^\top ~ V_{\rmf,k}^\top ~
Y_{\rmp,k}^\top ~ Y_{\rmf,k}^\top]^\top
$
where $
V_{\rmp,k} \in \R^{m\Tini \times q_k},~
V_{\rmf,k} \in \R^{m\Tf \times q_k},~
Y_{\rmp,k} \in \R^{n\Tini \times q_k},~
Y_{\rmf,k} \in \R^{n\Tf \times q_k}.
$
For each $k$, we get the prediction
$y_{[t, t+\Tf-1]}^\text{pred}(\hat V_{1:k}^{t})$ defined in \eqref{eq:subspace-predictor}, and keep only its first component $y_{t}^\text{pred}(\hat V_{1:k}^{t})$, simplifying its notation as $\hat y_t(k) \in \R^{n}$. 
Thus, by leveraging the nestedness of a flag, at each time $t$, we obtain $d$ different predictions $\hat{y}_t(1), \dots, \hat{y}_t(d)$ from matrices of different ranks based on only one matrix $\hat{V}_{1:d}^t$. 
Linking back to the motivation made in Section \ref{subsec:why_flag}, 
we compare against separately running \ac{FRONT} on $\Gr(16, q_k)$ for $k\in [1,d]$ and use the individually learned subspace $\hat{\bmU}_k^t \in \Gr(16, q_k)$ for prediction. In contrast, we use each $\hat{\bmV}_k^t$ with $k \in [1, d]$ from a single flag $\hat{\bmV}_{1:d}^t$ for prediction, requiring only one optimization instead of $d$ separate runs.
We measure the prediction performance by the cumulative prediction error over the whole simulation horizon $\Tsim-\Tf$ given by 
$\sum_{t=0}^{\Tsim-\Tf-1} \|y_t - \hat{y}_t\|_2^2,$
where $\hat{y}_t$ is the prediction at time $t$ from the corresponding model. We set $\Tini = \Tf = 4$ and the window size $T=20$. 

\begin{wrapfigure}[13]{r}{0.45\linewidth} 
\vspace{-20pt}
\centering
\begin{tikzpicture}[scale=0.65]
\begin{axis}[
    width=1.6\linewidth,
    height=0.95\linewidth,
    ybar,
    bar width=12pt,
    xlabel={Subspace Dimension},
    xtick=data,
    ymin=0,
    legend style={
        draw=none,
        at={(0.5, 1.05)},
        anchor=south,
        legend columns=-1,
        font=\large,
        column sep=6pt
    },
    legend image code/.code={
      \draw[draw=none, fill=#1] (0cm,-0.08cm) rectangle (0.35cm,0.18cm);
    },
    tick label style={font=\large},
    label style={font=\large},
    grid=both
]
\addplot[
    fill=orange!60, draw=none,
    error bars/.cd,
    y dir=both,
    y explicit,
]
table[
    x index=0, 
    y index=1,
    y error expr=\thisrowno{2}/2,
    col sep=comma, header=false,
]{data_prediction_experiment/flag_vs_gr_median.csv};
\addlegendentry{Nested subspaces in flag};

\addplot+[
    fill=myteal!40, draw=none,
    error bars/.cd,
    y dir=both,
    y explicit,
]
table[
    x index=0, 
    y index=3,
    y error expr=\thisrowno{4}/2,
    col sep=comma, header=false
]{data_prediction_experiment/flag_vs_gr_median.csv};
\addlegendentry{Grassmannians};
\end{axis}
\end{tikzpicture}
\captionsetup{width=0.85\linewidth}
\caption{Median cumulative prediction errors ($y$-axis): a single flag evaluated at different ranks (orange) versus individually learned Grassmann subspaces (green).
The whiskers span $30^{\text{th}}$ to $70^{\text{th}}$ percentiles.}
\label{fig:flag_vs_gr}
\end{wrapfigure}
Figure ~\ref{fig:flag_vs_gr} shows the median cumulative prediction errors of both approaches, where we set $d=8$ and $(q_1,\dots,q_d) = (8,\dots,15)$. The \ac{NSR} is set to be $\sigma=0.02$. 
All models were evaluated over 100 independent trials, and in each trial, all the models use the same offline input-output data $(u^\rmd_{[0,T_\rmd-1]}, y^\rmd_{[0,T_\rmd-1]} + \eta_{[0,T_\rmd]}^\rmd)$, input sequence $u_{[0, \Tsim-1]}$ and measurement noise $\eta_{[0,\Tsim-1]}$. 
The results demonstrate that the nested subspace predictions from a single flag achieve median prediction errors comparable to those of individually optimized Grassmannians across many dimensions.  We observe that the choice of subspace dimension can significantly affect performance. Optimal choice of subspace dimension is an area of future work. In the next part, we propose an ensembling strategy for mitigating the risk of poor subspace dimension selection.

\paragraph{Ensemble of Nested Subspace Prediction.}
Since the true system order is unknown, 
relying on a single subspace from the flag estimate may lead to biased predictions. To mitigate this, inspired by ensembling methods~\citep{dietterich2000ensemble}, we aggregate information from all candidate subspaces $\{\hat{\bmV}_{k}^t\}_{k=1}^d$. Specifically, we compute an ensemble of predictions by averaging $\hat{y}_t(k)$ obtained from \eqref{eq:subspace-predictor} across $k=1,\dots,d$, yielding
$
\hat{y}_t
 \coloneqq  
\frac{1}{d}\sum_{k=1}^d \hat{y}_t(k).
$
We use uniform averaging as a simple baseline ensemble rule. More sophisticated weighting schemes, for example based on the spectrum of the data covariance, may further improve performance and are an interesting direction for future work.
This approach leverages the hierarchical structure of the flag to balance the contribution of lower- and higher-dimensional subspaces, yielding more robust online predictions. We name this model as \textbf{average-$(q_1,\dots,q_d)$ prediction.}
We choose an ensemble of dimension $q_1 = 9$ and $q_2 = 10$ based on validation and evaluate the average-$(9, 10)$ prediction against several baseline models: no learning, N4SID, \ac{PAST}.
The \emph{no learning} model predicts future outputs without any adaptive updates to the subspace predictor.  
The N4SID method~\citep{van1994n4sid} is implemented in an online fashion: at each time $t$ we re-identify a state-space model $(A_t, B_t, C_t, D_t)$ from the current input-output data window $W_t=[w_{t-T+1} \cdots w_t]$. The state is estimated by a Kalman filter.  
Then the estimated state-space model is used to compute $\hat{y}_t$. 
\begin{figure}[tbp]
\centering
\pgfplotsset{compat=1.18}
\newcommand{\trajplot}[4]{%
  \begin{minipage}{0.48\linewidth}
  \centering
  \begin{tikzpicture}
  \begin{axis}[
    width=\linewidth, height=3.5cm,
    xmin=0, xmax=300,
    ymin=-8, ymax=8,
    grid=both,
    title={\footnotesize #1},
    tick label style={font=\scriptsize},
    label style={font=\scriptsize},
    legend style={at={(0.5,0.95)}, anchor=south, legend columns=-1,
      draw=none, fill=none, font=\scriptsize},
  ]
    \addplot[thick, color=gray]
      table[col sep=comma, header=false, x expr=\coordindex, y index=0]
      {#2};
    \addplot[thick, color=#4]
      table[col sep=comma, header=false, x expr=\coordindex, y index=0]
      {#3};
    \legend{True, Prediction}
  \end{axis}
  \end{tikzpicture}
  \end{minipage}%
}

\trajplot{No learning}
  {\detokenize{data_prediction_experiment/y_true_at_No_learning.csv}}
  {\detokenize{data_prediction_experiment/y_pred_No_learning.csv}}
  {myteal}
\hfill
\trajplot{$\mathrm{Flag}(16,(9,10))$}
  {data_prediction_experiment/y_true_at_Flag_9_10.csv}
  {data_prediction_experiment/y_pred_Flag_9_10.csv}
  {orange}

\caption{True and predicted trajectories.
Vertical axis: output. Horizontal axis: time steps.
Left: no learning baseline. Right: the average-$(9, 10)$ prediction model adapts and tracks the system after $t=\Tswitch$ quickly. }
\label{fig:y_prediction_two}
\end{figure}
The \ac{PAST} algorithm~\citep{yang1995projectionPAST} performs recursive subspace tracking based on each new $w_t$, where the estimated subspace is used in \eqref{eq:subspace-predictor} for prediction.
In addition, we evaluate the average-$(8,\dots,11)$ prediction to examine the performance when combining subspaces with relatively high cumulative prediction errors (see Figure \ref{fig:flag_vs_gr}). 
We vary the \ac{NSR} level as $\sigma \in \{0.01, 0.02, 0.05, 0.1\}$.
For each \ac{NSR} level, we conduct 100 independent trials.
All models are evaluated using the same 100 offline input-output data sets, identical input sequences, and the same noise realizations across all \ac{NSR} levels to ensure a fair comparison.
\paragraph{Results.}
\begin{wrapfigure}[13]{r}{0.5\linewidth} 
\centering
\pgfplotsset{compat=1.18}
\begin{tikzpicture}[scale=0.8]
\begin{axis}[
    width=1\linewidth, height=5cm,
    xmode=log, log basis x=10,
    xtick = {0.01, 0.02, 0.05, 0.1},
    xticklabels={0.01,0.02,0.05,0.1},
    xmin=0.01, xmax=0.1,
    ymin=5, ymax=11,
    ytick = {5, 10, 15, 20},
    grid=both, 
    legend columns=2,
    legend style={
        anchor=south,
        at={(0.5, 1.02)}, 
        draw=none,     
        font=\footnotesize,
    },
    tick label style={font=\footnotesize},
    label style={font=\footnotesize},
]

\pgfplotstableread[
  col sep=comma,
  header=false,
  skip first n=1,   
  trim cells=true
]{data_prediction_experiment/nsr_cpe_med_pct.csv}\datatable

\newcommand{\plotbandidx}[4]{%
  \addplot[name path=up, draw=none, forget plot]   
  table[x index=#1, y index=#3]{\datatable};
  \addplot[name path=down, draw=none, forget plot] 
  table[x index=#1, y index=#2]{\datatable};
  \addplot[fill opacity=0.15, draw=none, forget plot] 
  fill between[of=up and down];
  \addplot[fill=#4, fill opacity=0.2, draw=none, forget plot] 
  fill between[of=up and down];
}
\newcommand{\plotmedianidx}[5]{%
    \addplot+[
    line width=1.2pt, 
    mark=*, mark size=1.7pt, 
    color=#3,
    mark options={solid, fill=#3, draw=#3},
    #5
    ]
    table[x index=#1, y index=#2]{\datatable};
    \addlegendentry{#4}
}

\colorlet{N4SID}{gray}
\colorlet{PAST}{mycyan}
\colorlet{Gr8}{slategreen}
\colorlet{Gr9}{slategreen}
\colorlet{Gr10}{myblue}
\colorlet{Gr11}{seateal}
\colorlet{Gr15}{slategreen}
\colorlet{FlagMid}{myorange!80}
\colorlet{FlagBig}{myorange!80}
\colorlet{Flag910}{dustyrose}

\plotmedianidx{0}{4}{N4SID}{N4SID}{mark=square*}

\plotmedianidx{0}{7}{PAST}{PAST}{mark=triangle*}

\plotmedianidx{0}{13}{Gr8}{$\Gr(8)$}{solid}

\plotmedianidx{0}{16}{Gr9}{$\Gr(9)$}{mark=diamond}

\plotmedianidx{0}{19}{Gr10}{$\Gr(10)$}{mark=diamond*}

\plotmedianidx{0}{10}{Flag910}{$\Flag(16, (9,10))$}{mark=square}

\plotmedianidx{0}{22}{Gr11}{$\Gr(11)$}{mark=diamond*}

\plotmedianidx{0}{28}{FlagMid}{$\Flag(16, (8,\dots,11))$}{mark=square}

\end{axis}
\end{tikzpicture}
\captionsetup{width=0.8\linewidth}
\caption{Median cumulative prediction error ($y$-axis) versus varying \ac{NSR} ($x$-axis).}
\label{fig:NSR-CPE}
\end{wrapfigure}
Figure \ref{fig:y_prediction_two} shows, for the no learning baseline and the average-$(9, 10)$ prediction model, 
the predicted trajectory and the corresponding ground truth from the trial whose cumulative prediction error is closest to the median over 100 trials.
Figure \ref{fig:NSR-CPE} reports the median cumulative prediction error versus the \ac{NSR} for all models. The average-$(9, 10)$ prediction model achieves nearly the same median error as the  Grassmann models $\Gr(16, 9)$ and $\Gr(16, 10)$,
while requiring only a single optimization over one flag rather than separate optimizations for each dimension.
Since the range $(9, 10)$ is chosen a posteriori, this experiment mainly illustrates the computational advantage of a single flag optimization.
At all \ac{NSR} levels, the average-$(8, \dots, 11)$ model yields larger errors than the average-$(9, 10)$ model.
Since the dimensions in the average-$(8,\dots,11)$ model are chosen from a broader range without assuming knowledge of the true system order, this experiment reveals a limitation of the proposed method: including dimensions that underfit or overfit the system can degrade prediction performance. The results highlight a trade-off between robustness to unknown model order and prediction accuracy, as flags reduce the need for dimension selection but may degrade when the chosen dimension range is not well aligned with the true system complexity.
\vspace{-8pt}
\section{Conclusion}\label{sec:conclusion}
We proposed the \ac{FRONT} algorithm for online system identification and data-driven prediction that adapts to time-varying systems by recursively updating a nested hierarchy of subspaces. This enables the predictor to accommodate unknown or changing system dimensions without prior model information. 
We proved that GREAT is recovered as a special case of FRONT.
Numerical studies demonstrated that the proposed predictor effectively tracks abrupt changes in dynamics and consistently achieves better prediction error than most baselines even under high \ac{NSR}. 
Future work focuses on finite sample convergence of \ac{FRONT} for nontrivial flag signatures and integrating the proposed flag-based learning framework into an adaptive data-to-control pipeline.

\bibliography{l4dc2026-sample}

\end{document}